**Michiel Hazewinkel**     1     CWI
Direct line: +31-20-5924204     POBox 94079
Secretary: +31-20-5924233     1090GB Amsterdam
Fax: +31-20-5924166
E-mail: mich@cwi.nl

# The primitives of the Hopf algebra of noncommutative symmetric functions

by
*Michiel Hazewinkel*
*CWI*
*POBox 94079*
*1090GB Amsterdam*
*The Netherlands*

**Abstract**. Let *NSymm* be the Hopf algebra of noncommutative symmetric functions over the integers. In this paper a description is given of its Lie algebra of primitives over the integers, Prim(*NSymm*), in terms of recursion formulas. For each of the primitives of a basis of Prim(*NSymm*), indexed by Lyndon words, there is a recursively given divided power series over it. This gives another proof of the theorem that the algebra of quasi-symmetric functions is free over the integers.

**MSCS**: 16W30, 05E05, 17A50

**Key words and key phrases**: noncommutative symmetric functions, quasisymmetric functions, Lyndon word, Hopf algebra, primitive in a Hopf algebra, curve in a Hopf algebra, divided power series, Ditters-Shay bi-isobaric decomposition, Newton primitive, symmetric functions, Leibniz Hopf algebra, Lie Hopf algebra, free Lie algebra, graded Hopf algebra, coalgebra, free coalgebra, graded coalgebra, free associative algebra, Newton primitive, Verschiebung morphism, Frobenius morphism.

**1. Introduction**.
Let *NSymm* be the Hopf algebra of noncommutative symmetric functions, also known as the Leibniz Hopf algebra. As an algebra (over the integers) *NSymm* is simply the free algebra in countably many indeterminates:

$$NSymm = \mathbf{Z}\langle Z_1, Z_2, \cdots \rangle \qquad (1.1)$$

and the comultiplication is given by

$$\mu(Z_n) = \sum_{i+j=n} Z_i \otimes Z_j, \quad Z_0 = 1, \quad i, j \in \mathbf{N} \cup \{0\} \qquad (1.2)$$

*NSymm* is the noncommutative analogue of the Hopf algebra of symmetric functions

$$Symm = \mathbf{Z}[c_1, c_2 \cdots], \quad \mu(c_n) = \sum_{i+j=n} c_i \otimes c_j \qquad (1.3)$$

and more or less recently it has been discovered that very many of the remarkable structures and properties of the symmetric functions have natural noncommutative analogues in *NSymm* (or noncocommutative analogues in the graded dual *QSymm* of *NSymm*, the Hopf algebra of quasisymmetric functions); for instance, Schur functions, Newton primitives, representation theoretic interpretations, Frobenius reciprocity, ... ; see [3, 7, 8, 9, 10, 14],[6], and other papers. As often happens a number of things even become nicer or more transparent in the natural noncommutative generalization.

A primitive element in a Hopf algebra is an element $P$ such that

$$\mu(P) = 1 \otimes P + P \otimes 1 \qquad (1.4)$$

The set of primitives is a Lie algebra under the commutator difference product. For any Hopf algebra, and in particular for *NSymm*, there is interest in having a good description of its Lie algebra of primitives. One reason is that over a field of characteristic zero a cocommutative Hopf algebra is (isomorphic to) the universal enveloping algebra of its Lie algebra of primitives.

An outstanding question about *NSymm* is a good description of its Lie algebra of primitives over the integers, in particular writing down an explict basis of it as a free Abelian group. This is easy in the case of *Symm*, where the Lie algebra of primitives is commutative and a (rather canonical) basis is given by the Newton primitives (given by the same formula (1.7) below with the $Z_i$ replaced by the $c_i$). This is an instance of where the noncommutative version is more transparent and easier to prove than its commutative version. Formula (1.7) is clear and fits well with the recursion fomula

$$P_n(Z) = nZ_n - Z_{n-1}P_1(Z) - \cdots - Z_1P_{n-1}(Z)$$

In the commutative case the same recursion formula holds but it takes more than casual inspiration to guess that the coeffient of a monomial in the *c*'s in $P_n(c)$ is in fact the sum of the last indices of all noncommutative monomials that give rise to the same commutative monomial; and even so the explanation runs via noncommutative monomials.

This matter of primitives in the case of *NSymm* is vastly more complicated (and more interesting). To give some indications, lets first consider the matter over the rational numbers (which simplifies things quite a good deal). To do this consider yet another Hopf algebra

$$\mathcal{U} = \mathbf{Z}\langle U_1, U_2, \cdots \rangle, \quad \mu(U_n) = 1 \otimes U_n + U_n \otimes 1 \qquad (1.5)$$

In this case, the $U_n$ are primitives and the Lie algebra of primitives, Prim($\mathcal{U}$), is the free Lie algebra (over the integers) generated by the $U_n$. This is an object of sufficient interest and complexity that a book and more can be and has been devoted to it, [11].

Now, over the rationals *NSymm* and $\mathcal{U}$ are isomorphic as Hopf algebras. The isomorphism is given by setting

$$1 + Z_1 t + Z_2 t^2 + Z_3 t^3 + \cdots = \exp(U_1 t + U_2 t^2 + U_3 t^3 + \cdots) \qquad (1.6)$$

which gives formulas for the $Z_n$ as polynomials in $U_1, U_2, \cdots, U_n$, which are then used to define an algebra morphism of *NSymm* into $\mathcal{U}$, which turns out to be a Hopf algebra isomorphism; see [4] for two proofs of the latter fact.

Thus, up to isomorphism, there is a good description of the primitives of *NSymm* over the rationals. Actually, one can do even better (from certain points of view). Define the (noncommutative) Newton primitives by

$$P_n(Z) = \sum_{r_1 + \cdots + r_k = n} (-1)^{k+1} r_k Z_{r_1} Z_{r_2} \cdots Z_{r_k}, \quad r_i \in \mathbf{N} = \{1, 2, \cdots\} \qquad (1.7)$$

It is easily proved by induction that the $P_n(Z)$ are primitives of *NSymm*, and it is also easy to see that over the rationals *NSymm* is the free associative algebra generated by the $P_n(Z)$. Thus over the rationals the Lie algebra of primitives of *NSymm* is simply the free Lie algebra generated by the $P_n(Z)$. Over the integers things are vastly different. For one thing Prim(*NSymm*) is most definitely not a free Lie algebra; rather it tries to be something like a divided power Lie algebra (though I do not know what such a thing would be). As it turns out there is a basis of Prim(*NSymm*) that includes the $P_n(Z)$ and these generate a free sub Lie algebra of Prim(*NSymm*). This free sub Lie algebra is a rather small part of Prim(*NSymm*). Both are infinite dimensional free Abelian groups, and the one is a full rank subgroup of the other (meaning that over the rationals they span the same vectorspace). So, to indicate how small one is within the other we have to consider some of the extra structure that is present, viz that we are dealing with *graded* Abelian groups. The grading is given by weight. In *NSymm* give $Z_n$ weight $n$, and a monomial $Z_{a_1} Z_{a_2} \cdots Z_{a_m}$ weight $a_1 + a_2 + \cdots + a_m$. This makes *NSymm* a graded Abelian group. If an element $P$ is primitive so are all its homogeneous components. Thus Prim(*Symm*) is also a graded Abelain group. Further $P_n(Z)$, as defined by formula (1.7) above, is homogeneous of weight $n$. Thus the free Lie algebra FL($P$) generated by them is also graded. The homogenous components of both Prim(*NSymm*) and FL($P$) are free of finite rank, FL($P$)$_n \subset$ Prim(NSymm)$_n$ and we can study the value of the index of one in the other. The results for the first few values of $n$ are as follows

| n     | 1 | 2 | 3 | 4 | 5   | 6     |   |
|-------|---|---|---|---|-----|-------|---|
| Index | 1 | 1 | 2 | 6 | 576 | 69120 |   |

,
a sequence of numbers that grows far faster than any exponential. There is in fact a formula. For any word $\alpha = [a_1, a_2, \cdots, a_n]$ over the natural numbers let $g(\alpha)$ be the gcd of its entries and let $k(\alpha)$ be the product of its entries, then

$$\text{Index of FL}(P)_n \text{ in Prim}(NSymm)_n = \sum_{\alpha \in LYN, \, \text{wt}(\alpha)=n} \frac{k(\alpha)}{g(\alpha)}.$$

Here *LYN* is the set of Lyndon words over **N** (defined below).

To describe the basis of Prim(*NSymm*) alluded to (of Prim(*NSymm*) as a free Abelian group) a number of definitions are needed.

A *word* over the natural numbers **N** = {1,2,$\cdots$} is simply a sequence of natural numbers $\alpha = [a_1, a_2, \cdots, a_m]$. The length, lg($\alpha$), of the word $\alpha$ is $m$, and its weight is wt($\alpha$) = $a_1 + a_2 + \cdots + a_m$. A word over **N** of weight $n$ is also a called a *composition* of $n$. The proper tails of $\alpha$ are the words $[a_i, a_{i+1}, \cdots, a_m]$, $i = 2, \cdots, m$. A word $\alpha = [a_1, \cdots, a_m]$ is lexicographically larger than a word $\beta = [b_1, \cdots, b_n]$ iff there is an index $i$ ($1 \leq i \leq \min\{m,n\}$) such that $a_1 = b_1$, $\cdots$, $a_{i-1} = b_{i-1}$, $a_i > b_i$, or $m > n$ and $a_j = b_j$, $j = 1, \cdots, n$. A word is *Lyndon* if it is lexicographicaly smaller than each of its proper tails.

Let $H$ be a Hopf algebra. A *curve* in $H$, also called a *DPS* (*divided power series*), is a sequence of elements from $H$, $d = (d(0) = 1, d(1), d(2), \cdots)$ such that for all $n$

$$\mu_H(d(n)) = \sum_{i+j=n} d(i) \otimes d(j) \tag{1.8}$$

Here $\mu_H$ is the comultiplication of the Hopf algebra $H$. An often convenient notation for a curve is $d(t) = 1 + d(1)t + d(2)t^2 + \cdots \in H[[t]]$; i.e. it is written as a power series in $t$ over $H$ (where $t$ is a (counting) variable commuting with all elements of $H$).
With this notation, if $d(t)$ is a curve than so is $d(t^r)$, which corresponds to the sequence

$$1,\underbrace{0,\cdots,0}_{r-1},d(1),\underbrace{0,\cdots,0}_{r-1},d(2),\cdots.$$

And if $d(t), d_1(t), d_2(t)$ are curves than so are the power series product $d_1(t)d_2(t)$ and the power series inverse $d(t)^{-1}$. The latter corresponds to the sequence $1, \iota(d(0)), \iota(d(2)), \cdots$ where $\iota$ is the antipode of $H$.

Now consider the free product *2NSymm* of *NSymm* with itself

$$2NSymm = \mathbf{Z}\langle X_1, Y_1, X_2, Y_2, \cdots \rangle, \quad \mu(X_n) = \sum_{i+j=n} X_i \otimes X_j, \quad \mu(Y_n) = \sum_{i+j=n} Y_i \otimes Y_j \quad (1.9)$$

In *2NSymm* there are two obvious (natural) curves, viz

$$X(s) = 1 + X_1 s + X_2 s^2 + \cdots, \quad Y(t) = 1 + Y_1 t + Y_2 t^2 + \cdots \quad (1.10)$$

(where $s$ is a second counting variable). Now, consider the commutator product $X(s)^{-1}Y(t)^{-1}X(s)Y(t)$. The Shay-Ditters bi-isobaric decomposition theorem says that this commutator product can be written uniquely as an ordered product

$$X(s)^{-1}Y(t)^{-1}X(s)Y(t) = \prod_{\gcd(a,b)=1} (1 + L_{a,b}(X,Y)s^a t^b + L_{2a,2b}(X,Y)s^{2a}t^{2b} + \cdots) \quad (1.11)$$

Here the ordered product is over all pairs of natural numbers $(a,b) \in \mathbf{N} \times \mathbf{N}$ with greatest common divisor, gcd, equal to 1, and the ordering is

$$(a,b) > (a',b') \iff a+b > a'+b' \text{ or } (a+b = a'+b' \text{ and } a > a') \quad (1.12)$$

Actually it does not matter much which ordering is used. This decomposition theorem is pretty obvious once one observes that the monomials of the form $s^a$ and $t^b$ in the commutator product have coefficient zero and that for each $(k,l) \in \mathbf{N} \times \mathbf{N}$ there is precisely one pair $(a,b), \gcd(a,b) = 1$, such that the monomial $s^k t^l$ occurs in

$$d_{a,b}(s,t) = 1 + L_{a,b} s^a t^b + L_{2a,2b} s^{2a} t^{2b} + \cdots \quad (1.13)$$

viz $(a,b) = (k/\gcd(k,l), l/\gcd(k,l))$. The $L_{k,l}(X,Y)$ have weight $k$ in $X$ and weight $l$ in $Y$.

Writing (1.11) in the form

$$X(s)Y(t) = Y(t)X(s) \prod_{\gcd(a,b)=1} (1 + L_{a,b}(X,Y)s^a t^b + L_{2a,2b}(X,Y)s^{2a}t^{2b} + \cdots)$$

there immediately follows a recursion formula for the $L_{ra,rb}(X,Y)$; see below in section 3 for a description of the explicit formula.

It follows readily from (1.11) that for all $(a,b) \in \mathbf{N} \times \mathbf{N}$, $\gcd(a,b) = 1$, the sequences

$$1, L_{a,b}(X,Y), L_{2a,2b}(X,Y), \cdots \quad (1.14)$$

are curves in *2NSymm*. It follows also, because of the nature of the comultiplication (1.9) of *2NSymm*, that the curves (1.14) can be used to define a new curve from two old ones: if $d_1, d_2$ are two curves in a Hopf algebra $H$, then for each $(a,b) \in \mathbf{N} \times \mathbf{N}, \gcd(a,b) = 1,$

$$d_{a,b}(d_1,d_2) = (1, L_{a,b}(d_1,d_2), L_{2a,2b}(d_1,d_2), \cdots) \qquad (1.15)$$

is also a curve in $H$.

There is now sufficient notation to give an explicit recursive description and construction of a basis for Prim(*NSymm*). The basis consists of primitives $P_\alpha$, one for each Lyndon word $\alpha$. To each Lyndon word there are associated three things, viz a number $g(\alpha)$, the gcd of the entries of $\alpha = [a_1, \cdots, a_m]$, a curve $d_\alpha$, and a primitive $P_\alpha$ that is homogeneous of weight wt($\alpha$). The curves associated to $\alpha$ and $r\alpha = [ra_1, \cdots, ra_m]$ are the same. Here is the recursive description. For $\alpha$ of length 1, i.e. $\alpha = [n]$, we have

$$g([n]) = n, \quad d_{[n]} = (1, Z_1, Z_2, \cdots), \quad P_{[n]} = P_n(Z) \qquad (1.16)$$

For a Lyndon word $\alpha$ of length $>1$ let $\alpha'' = [a_i, \cdots, a_m]$ be its lexicographically smallest proper tail (suffix), and let $\alpha' = [a_1, \cdots, a_{i-1}]$ be the corresponding prefix of $\alpha$. Then both $\alpha', \alpha''$ are Lyndon words. Of course $g(\alpha) = \gcd(g(\alpha'), g(\alpha''))$. The entities associated to $\alpha$ are now

$$\begin{aligned} g(\alpha) &= \gcd(g(\alpha'), g(\alpha'')) = \gcd(a_1, \cdots a_m), \\ d_\alpha &= (1, L_{g(\alpha')/g(\alpha), g(\alpha'')/g(\alpha)}(d_{\alpha'}, d_{\alpha''}), L_{2g(\alpha')/g(\alpha), 2g(\alpha'')/g(\alpha)}(d_{\alpha'}, d_{\alpha''}), \cdots), \\ P_\alpha &= P_{g(\alpha)}(d_\alpha) \end{aligned} \qquad (1.17)$$

The main theorem is

**Theorem**. The $P_\alpha$, where $\alpha$ runs over all Lyndon words, form a basis (over the integers) of Prim(*NSymm*).

From the construction above it is immediate that for Lyndon words $\alpha$ with $g(\alpha) = 1$, the corresponding primitive $P_\alpha$ is the first term of a curve (DPS). If $g(\alpha) > 1$, there is also such a curve. To see that there is a second bi-isobaric decomposition theorem. This time we work in *NSymm* itself and consider

$$Z(s)^{-1} Z(t)^{-1} Z(s+t) \quad NSymm[[s,t]], \quad Z(s) = 1 + Z_1 s + Z_2 s^2 + \cdots \qquad (1.18)$$

Again it is clear that the monomials $s^k$ and $t^l$ have coefficient zero in $Z(s)^{-1} Z(t)^{-1} Z(s+t)$, and again there is a natural bi-isobaric decomposition

$$Z(s)^{-1} Z(t)^{-1} Z(s+t) = \prod_{\gcd(a,b)=1} (1 + N_{a,b}(Z) s^a t^b + N_{2a,2b}(Z) s^{2a} t^{2b} + \cdots) \qquad (1.19)$$

And again it follows that for each pair $(a,b) \in \mathbf{N} \times \mathbf{N}$, $\gcd(a,b) = 1$ the sequences

$$1, N_{a,b}(Z), N_{2a,2b}(Z), \cdots \qquad (1.20)$$

are curves, and again these curves can be used to construct a new one from a known one. In particular, if $d$ is a curve in *NSymm*, then for each $n > 1$

$$1, N_{1,n-1}(d), N_{2,2n-2}(d), \cdots \qquad (1.21)$$

is again a curve. A quick check shows that

$$N_{1,n-1}(Z) = P_n(Z) \tag{1.22}$$

as defined by (1.7) above. Thus for all Lyndon words $\alpha$ the associated primitive $P_\alpha$ is the first term of a curve (DPS). It follows that the graded dual Hopf algebra of *NSymm*, the algebra *QSymm* of quasisymmetric functions, is free as an algebra over the integers. A proof of the implication (all primitives in *NSymm* extend to curves (DPS's)) $\Rightarrow$ *QSymm* is free) is given in an appendix (for the case at hand. It is adapted from the one in [12] and the one in [13] (for the case of Hopf algebras over a field of characteristic zero) and uses ingredients from both. So

**Theorem**. The Hopf algebra of quasisymmetric functions, the graded dual of *NSymm*, as an algebra is (commutative) free.

For another proof of this theorem and more about the algebra *QSymm* of quasisymmetric functions as the dual of *NSymm*, see [5, 6].

**Acknowledgements**. The description of the primitives of *NSymm* given in this paper is essentially the same as the one in the preprint [12]. There are sign differences and the set of words used by Brian Shay is quite different from the set of Lyndon words used here. So the actual explicit formulas, when written out completely, are quite a bit different. These differences probably do not really matter. Any Hall-like or Lazard-like set of words should work. Also different orderings of the set $\{(r,s): r,s \in \mathbf{N} \setminus \{0\}, \gcd(r,s) = 1\}$ can no doubt be used.
     One of the most essential ingredients of the construction, the Ditters-Shay bi-isobaric decomposition theorem, cf above and below, is due, independently, to both Shay and Ditters; see, [1, 2, 12] where it occurs in somewhat different forms than here. Another nice notion, not crucial but very nice and useful, that of a **V**-curve, see below, is also due to Ditters.
     The preprint [12] is very difficult to decipher; first because of the horrendous notations used and second because of a dozen or more typos and/or inaccuracies on practically all the more important pages. Still, currently, I have the impression that the preprint is basically correct. In that case the first proof that the algebra of quasisymmetric functions is free as an algebra over the integers is due to Shay. The second proof is in [5] and the present paper provides a third one based on rather similar ideas as those in [12], but with quite different proofs and some new constructions such as the second bi-isobaric decomposition theorem.

## 2. Curves, 2-curves, and V-curves.

2.1. *Curves.*
     A curve in a Hopf algebra $H$ is a sequence of elements

$$d = (d(0) = 1, d(1), d(2), \cdots) \tag{2.1.1}$$

such that

$$\mu_H(d(n)) = \sum_{i+j=n} d(i) \otimes d(j) \tag{2.1.2}$$

A curve in a Hopf algebra is also called a divided power series. When a curve or DPS is written as a sequence like (2.1.1), the term $d(0)=1$ is often omitted. Note that $d(1)$ is a primitive.
     A convenient way, in many situations, to write a curve is as a power series in a counting variable $t$ (which commutes with all elements of $H$)

$$d(t) = 1 + d(1)t + d(2)t^2 + \cdots \tag{2.1.3}$$

Noncommutative multiplication and inversion of power series with coefficients in $H$ turns the set of curves in a Hopf algebra into a (usually noncommutative) group. In the case of commutative formal groups these groups, enriched with a number of functorial operations on them, are classifying. It is a matter of absorbing interest to investigate to what extent this may still be true in more general situations. The inverse of $d(t)$ in (2.1.3) is the power series

$$d^{-1}(t) = 1 + \iota_H(d(1))t + \iota_H(d(2))t^2 + \cdots \tag{2.1.4}$$

where $\iota_H$ is the antipode of the Hopf algebra $H$.

There are two more useful ways of looking at curves. Let $\mathcal{C}$ be the coalgebra

$$\mathcal{C} = \mathbf{Z} \bigoplus_{i=1} \mathbf{Z}Z_i, \quad \mu(Z_n) = \sum_{i+j=} Z_i \otimes Z_j, \quad Z_0 = 1 \tag{2.1.5}$$

With the co-unit given by projection on the zero-th factor. Then a curve in a Hopf algebra $H$ is exactly the same as a coalgebra morphism

$$\mathcal{C} \longrightarrow H, \quad Z_i \mapsto d(i) \tag{2.1.6}$$

Further the multipication of two curves $d, d'$ with corresponding coalgebra morphisms $\varphi, \varphi'$ coresponds to the convolution

$$\mathcal{C} \xrightarrow{\mu} \mathcal{C} \otimes \mathcal{C} \xrightarrow{\varphi \otimes \varphi'} H \otimes H \xrightarrow{m} H \tag{2.1.7}$$

Finally, the Hopf algebra of noncommutative symmetric functions $NSymm = \mathbf{Z}\langle Z_1, Z_2, \cdots \rangle$ is the free associative algebra on the $Z_i$. Thus a curve in $H$ is also the same as a Hopf algebra morphism

$$NSymm \xrightarrow{\varphi} H, \quad \varphi(Z_i) = d(i) \tag{2.1.8}$$

and this is the point of view that shall be frequently used below. There is, however, some danger in this. It is very tempting to write down a similar diagram as (2.1.7)

$$NSymm \xrightarrow{\mu} NSymm \otimes NSymm \xrightarrow{\varphi \otimes \varphi'} H \otimes H \xrightarrow{m} H \tag{2.1.9}$$

and to think that this is the morphism of Hopf algebras corresponding to the product of $d$ and $d'$. As a rule it is not if $H$ is noncommutative; in particular this is not the case for $H = NSymm$. The problem is that if $H$ is noncommutative then $H \otimes H \xrightarrow{m} H$ is not an algebra morphism.

The terminology 'curve' comes from a special case. Let $F$ be a formal group over a ring A, $R(F) = A[[X_1, X_2, \cdots, X_n]]$ its contravariant bialgebra, and $U(F)$ its covariant Lie algebra. A curve in $U(F)$ is a coalgeba morphism $\mathcal{C} \longrightarrow U(F)$; duality gives an algebra morphism $R(F) \longrightarrow A[[t]]$, i.e. a curve in the sense of formal geometry.

2.2. *2-Curves.*
A 2-curve in a Hopf algebra $H$ is a collection of elements $c(n,m) \in H$, indexed by pairs of nonnegative integers $(n,m) \in \mathbf{N} \cup \{0\} \times \mathbf{N} \cup \{0\}$, $c(0,0) = 1$, such that

$$\mu(c(n,m)) = \sum_{\substack{n_1+n_2=n \\ m_1+m_2=m}} c(n_1, m_1) \otimes c(n_2, m_2) \tag{2.2.1}$$

A convenient way of writing a 2-curve is as a powerseries in two variables

$$c(s,t) = 1 + \sum_{n+m>0} c(n,m) s^n t^m \qquad (2.2.2)$$

Power series products and power series inverses of 2-curves are again 2-curves. As in the case of curves (i.e. 1-curves) there are interpretations in terms of coalgebra morphisms and Hopf algebra morphisms, but those will not be needed here.

A curve

$$d(s) = 1 + d(1)s + d(2)s^2 + \cdots \qquad (2.2.3)$$

can be seen as a (degenerate) 2-curve. An example of a most important 2-curve is

$$X(s)^{-1} Y(t)^{-1} X(s) Y(t) \in 2NSymm[[s,t]]$$

where

$$X(s) = 1 + X_1 s + X_2 s^2 + \cdots \quad \text{and} \quad Y(t) = 1 + Y_1 t + Y_2 t^2 + \cdots$$

are the two (canonical) natural curves in $2NSymm$, see (1.9) above.

2.2.4. **Lemma.** Let $a,b \in \mathbf{N} \cup \{0\}$, not both zero, and let $d(s,t)$ be a two curve of the form

$$d(s,t) = 1 + c_1 s^a t^b + c_2 s^{2a} t^{2b} + \cdots \qquad (2.2.5)$$

Then $c(t) = 1 + c_1 t + c_2 t^2 + \cdots$ is a curve.

Proof. Because $d(s,t)$ is a 2-curve we have

$$\mu(d(n,m)) = \sum_{\substack{n_1+n_2=n \\ m_1+m_2=m}} d(n_1, m_1) \otimes d(n_2, m_2) \qquad (2.2.6)$$

Now most of the terms on the right hand side of (2.2.6) are zero. The only ones that are possibly nonzero are of the form $(n_1, m_1) = r_1(a,b)$, $(n_2, m_2) = r_2(a,b)$ and then $(n,m) = (r_1+r_2)(a,b) = r(a,b)$, $r_1+r_2 = r$ and $d(n,m) = c_r$, $d(n_1,m_1) = c_{r_1}$, $d(n_2,m_2) = c_{r_2}$ so that

$$\mu(c_r) = \sum_{r_1+r_2=r} c_{r_1} \otimes c_{r_2}$$

proving that $c(t)$ is a curve.

2.3. **V**-*curves*.

On *NSymm* and $2NSymm$ there are some some remarkable Hopf algebra endomorphisms, called Verschiebung. There is one for every $r \in \mathbf{N}$, and they are defined as follows (on *NSymm* and $2NSymm$ respectively):

$$\mathbf{V}_r(Z_n) = \begin{cases} Z_{n/r} & \text{if } r \text{ divides } n \\ 0 & \text{otherwise} \end{cases} \qquad (2.3.1)$$

$$\mathbf{V}_r(X_n) = \begin{cases} X_{n/r} & \text{if } r \text{ divides } n \\ 0 & \text{otherwise} \end{cases}, \quad \mathbf{V}_r(Y_n) = \begin{cases} Y_{n/r} & \text{if } r \text{ divides } n \\ 0 & \text{otherwise} \end{cases} \quad (2.3.2)$$

A curve

$$d = (d(0), d(1), d(2), \cdots)$$

in *NSymm* or 2*NSymm* is a **V**-curve if

$$\mathbf{V}_r(d(n)) = \begin{cases} d(n/r) & \text{if } r \text{ divides } n \\ 0 & \text{otherwise} \end{cases}$$

Let $\varphi: NSymm \to NSymm$, respectively, $\varphi: NSymm \to 2NSymm$ be the morphism of Hopf algebras corresponding to the curve $d$ in *NSymm*, respectively, 2*NSymm*. Then $d$ is a **V**-curve if and only if $\varphi$ commutes with the endomorphisms $\mathbf{V}_r$.

2.4. *Substituting curves in curves.*

Let $d$ be a curve in *NSymm*, and $d'$ a curve in a Hopf algebra $H$. Define $d(d')$ as the sequence of elements of $H$ obtained by replacing the $Z_1, Z_2, \cdots$ in $d$ with $d'(1), d'(2), \cdots$. More precisely, $d(i)$ is a polynomial in the $Z_j$, $d(i)(Z_1, Z_2, \cdots)$, and then

$$d(d')(i) = d(i)(d'(1), d'(2), \cdots)$$

2.4.1. Proposition. As above let $d$ be a curve in *NSymm* and $d'$ a curve in a Hopf algebra $H$. Then
    (i) $d(d')$ is a curve in $H$.
    (ii) When $H = NSymm$ and both $d, d'$ are **V**-curves, then $d(d')$ is a **V**-curve.

Proof. Let $\varphi_d: NSymm \to NSymm$ and $\varphi_{d'}: NSymm \to H$ be the morphisms of Hopf algebras corresponding to $d, d'$. Then the morphism of Hopf algebras coresponding to $d(d')$ is the composed morphism

$$NSymm \xrightarrow{\varphi_d} NSymm \xrightarrow{\varphi_{d'}} H \quad (2.4.2)$$

Composing Hopf algebra morphisms gives a Hopf algebra morphism, so $d(d')$ is a curve. If $H = NSymm$ and $d, d'$ are **V**-curves, then $\varphi_d, \varphi_{d'}$ both commute with the $\mathbf{V}_r$, and hence so does their composed morphism (2.4.2), proving that $d(d')$ is a **V**-curve.

Now let $c$ be a curve in 2*NSymm*, and $d, d'$ curves in a Hopf algebra $H$. The sequence $c(d, d')$ is obtained from $c$ by replacing the $X_i$ in $c$ by $d(i)$ and the $Y_j$ by $d'(j)$.

2.4.3. Proposition. As above let $c$ be a curve in 2*NSymm*, and $d, d'$ curves in a Hopf algebra $H$. Then
    (i) $c(d, d')$ is a curve in $H$
    (ii) When $H = NSymm$ and all three curves $c, d, d'$ are **V**-curves, then $c(d, d')$ is a **V**-curve.

Proof. The pair of curves $d, d'$ defines a Hopf algebra morphism

$$\varphi_{d,d'}: 2NSymm \to H, \quad X_i \mapsto d_i, Y_j \mapsto d'_j \quad (2.4.4)$$

Then $c(d,d')$ corresponds to the composed morphism

$$NSymm \xrightarrow{\varphi_c} 2NSymm \xrightarrow{\varphi_{d,d'}} H \qquad (2.4.5)$$

and as a composition of Hopf algebra morphisms this is a Hopf algebra morphism, so that $c(d,d')$ is a curve. If $H = NSymm$ and $d, d'$ are **V**-curves, then $\varphi_{d,d'}$ commutes with the $\mathbf{V}_r$ (on $2NSymm$ and $NSymm$ respectively). Hence the composed morphism (2.4.5) also commutes with the $\mathbf{V}_r$ making $c(d,d')$ a **V**-curve.

   2.4.6. **Comments.** Write $\mathrm{Curve}(H)$ for the group of curves in a Hopf algebra $H$ where the group multiplication is multiplication of power series. Write $\mathcal{E} = \mathrm{End}_{Hopf}(NSymm)$, then the first construction above amounts to defining a right action of the semigroup $\mathcal{E}$ on $\mathrm{Curve}(H)$. For $d \in \mathrm{Curve}(H)$ and $\psi \in \mathcal{E}$, the curve $d\psi$ is the one corresponding to the composition of Hopf algebra morphisms

$$NSymm \xrightarrow{\psi} NSymm \xrightarrow{\varphi_d} H \qquad (2.4.7)$$

If $H$ is commutative this action respects the group structure on $\mathrm{Curve}(H)$. If $H$ is not commutative, in particular when $H = NSymm$, this is not the case. This here is an instance of the possible pitfalls in looking at curves as Hopf algebra morphisms. If, erroneously, one took (2.1.9) as corresponding to the power series product of curves, it would follow that this right action does respect the group structure.

   In case $H = NSymm$ there is also a left action. For $d \in \mathrm{Curve}(NSymm)$, $\psi \in \mathcal{E}$, $\psi d$ is the curve coresponding to the composition of Hopf algebra morphisms

$$NSymm \xrightarrow{\varphi_d} NSymm \xrightarrow{\psi} NSymm \qquad (2.4.8)$$

This one does respect the group structure. In this case of course

$$\mathrm{End}_{Hopf}(NSymm) = \mathcal{E} = \mathrm{Curve}(NSymm)$$

and $\mathrm{Curve}(NSymm)$ is a set with a noncommutative addition on it (power series multiplication of curves), a noncommutative multiplication on it (composition of endomorphisms) and the multiplication is distributive over the addition on the left but not on the right. There is also a unit (the identity endomorphism, or, as a curve, the natural curve $1 + Z_1 t + Z_2 t^2 + \cdots$).

## 3. Isobaric decomposition.

Consider again

$$2NSymm = \mathbf{Z}\langle X_1, Y_1, X_2, Y_2, \cdots \rangle, \quad \mu(X_n) = \sum_{i+j=n} X_i \otimes X_j, \quad \mu(Y_n) = \sum_{i+j=n} Y_i \otimes Y_j \qquad (3.1)$$

and the two natural curves

$$X(s) = 1 + X_1 s + X_2 s^2 + \cdots, \quad Y(t) = 1 + Y_1 t + Y_2 t^2 + \cdots \qquad (3.2)$$

and consider the commutator product

$$X(s)^{-1} Y(t)^{-1} X(s) Y(t) \qquad (3.3)$$

On the set of pairs of nonnegative integers consider the ordering

$$(u,v) <_{wl} (u',v') \iff u+v < u'+v' \text{ or } (u+v = u'+v' \text{ and } u < u') \quad (3.4)$$

(Here the index *wl* on $<_{wl}$ is supposed to be a mnemonic for weight first, then lexicographic.)

3.5. **Theorem** (Ditters-Shay bi-isobaric decomposition theorem). There are 'higher commutators' (or perhaps better 'corrected commutators')

$$L_{u,v}(X,Y) \in \mathbf{Z}\langle X,Y \rangle, \quad (u,v) \in \mathbf{N} \times \mathbf{N} \quad (3.6)$$

such that

$$X(s)^{-1}Y(t)^{-1}X(s)Y(t) = \prod_{\gcd(a,b)=1} (1 + L_{a,b}(X,Y)s^a t^b + L_{2a,2b}(X,Y)s^{2a}t^{2b} + \cdots) \quad (3.7)$$

where the product is an ordered product for the ordering $<_{wl}$ just introduced, (3.4). Moreover

(i) $L_{u,v}(X,Y) = [X_u, Y_v] + \text{(terms of length} \geq 3)$ \quad (3.8)

(ii) $L_{u,v}(X,Y)$ is homogeneous of weight $u$ in $X$ and of weight $v$ in $Y$. \quad (3.9)

(iii) For $\gcd(a,b) = 1$, $1 + L_{a,b}(X,Y)s^a t^b + L_{2a,2b}(X,Y)s^{2a}t^{2b} + \cdots$ is a curve. It is also a **V**-curve.

**Proof.** All this basically follows from two simple observations. First that putting $s$ or $t$ zero in the left hand side of (3.7) gives 1, so that there are no pure powers of $s$ or $t$ in (3.7); second that for each $(u,v) \in \mathbf{N} \times \mathbf{N}$ there is precisely one $(a,b) \in \mathbf{N} \times \mathbf{N}$, $\gcd(a,b) = 1$, such that $s^u t^v$ occurs in

$$1 + L_{a,b}(X,Y)s^a t^b + L_{2a,2b}(X,Y)s^{2a}t^{2b} + \cdots \quad (3.10)$$

viz $(a,b) = (u,v) / \gcd(u,v)$. In more detail, rewrite (3.7), as

$$X(s)Y(t) = Y(t)X(s) \prod_{\gcd(a,b)=1} (1 + L_{a,b}(X,Y)s^a t^b + L_{2a,2b}(X,Y)s^{2a}t^{2b} + \cdots) \quad (3.11)$$

Comparing coefficients of $s^u t^v$ left and right one finds

$$X_u Y_v = Y_v X_u + \sum_{\substack{u_0,v_0 \geq 0; k \geq 1 \\ u_i,v_i \geq 1, i=1,\cdots,k \\ u_0+u_1+\cdots+u_k=u \\ v_0+v_1+\cdots+v_k=v}} Y_{v_0} X_{u_0} \prod_i L_{u_i,v_i}(X,Y) \quad (3.12)$$

where the product is an ordered one for the ordering

$$(u_1,v_1)/\gcd(u_1,v_1) \leq_{wl} (u_2,v_2)/\gcd(u_2,v_2) \leq_{wl} \cdots \leq_{wl} (u_k,v_k)/\gcd(u_k,v_k) \quad (3.13)$$

This is really a recursion formula for $L_{u,v}(X,Y)$ (this term being the case $u_0 = 0, v_0 = 0, k = 1$). Explicitely

$$L_{u,v}(X,Y) = [X_u, Y_v] - \sum_{\substack{u_0, v_0 \geq 0; k \geq 1 \\ u_i, v_i \geq 1, i=1,\cdots,k \\ u_0+u_1+\cdots+u_k=u \\ v_0+v_1+\cdots+v_k=v \\ u_0+v_0 \geq 1 \text{ or } k \geq 2}} Y_{v_0} X_{u_0} \prod_i L_{u_i,v_i}(X,Y) \tag{3.14}$$

The last restriction in the sum in (3.14) is simply a way of saying that every term under the sum sign has at least two factors (and not of the form $Y_{v_0} X_{u_0}$ by the first condition in the sum of (3.14)) so that all the $L_{u',v'}(X,Y)$ on the right hand side of (3.14) have lower weight, $u'+v' < u+v$.

Now define the $L_{u,v}(X,Y)$ by formula (3.14), then (3.12) holds, and hence (3.11) and (3.7). This takes care of existence (and uniqueness for that matter). Statements (i) and (ii) follow immediately from the recursion formula (3.14).

To show that the series (3.10) are curves use induction with respect to bidegree for the wl-ordering. Suppose that that all the factors on the right of (3.7) have been shown to be 2-curves up to (but not including) bidegree $(u,v)$. Let $(a,b) = (u,v)/\gcd(u,v)$. For all the power series

$$1 + L_{a,b}(X,Y)s^a t^b + L_{2a,2b}(X,Y)s^{2a} t^{2b} + \cdots \tag{3.15}$$

with $(a',b') \neq (a,b)$, the coefficient of $s^u t^v$ is zero. Also a term from

$$\sum_{\substack{u'+u''=u \\ v'+v''=v}} \text{coeff}(s^{u'} t^{v'}) \cdot \text{coeff}(s^{u''} t^{v''}) \tag{3.16}$$

in (3.15) can be nonzero only if $(u',v')$ and $(u'',v'')$ are both multiples of $(a',b')$. But that would imply that $(u,v)$ is a multiple of $(a',b')$, which is not the case. Thus all the factors on the right of (3.7), except possibly (3.10), are 2-curves up to and including didegree $(u,v)$. Because the left hand side of (3.7) is a 2-curve, it follows that the last remaining term, (3.10), is also a 2-curve up to and including bidegree $(u,v)$.

The proof that the series (3.10) are **V**-curves goes exactly the same way. Again suppose that this has been proved up to bidegree $(u,v)$. The coefficient of $s^u t^v$ is zero and so is the coefficient of all $s^{u'} t^{v'}$ for $(u',v') = r^{-1}(u,v)$, $r \mid \gcd(u,v)$. Thus all the factors on the right of (3.7) are **V**-curves up to and including bidegree $(u,v)$ except possibly (3.10) itself. But the left hand side of (3.7) is a **V**-curve. Hence (3.10) is also a **V**-curve up to and including bidegree $(u,v)$.

Before stating the second isobaric decomposition theorem some preparation is needed. Consider the natural curve

$$Z(t) = 1 + Z_1 t + Z_2 t^2 + \cdots$$

in *NSymm*.

**3.17. Lemma.** $Z(s+t)$ is a 2-curve.

Proof. The coefficient of $s^a t^b$ in $Z(s+t)$ is $\binom{a+b}{a} Z_{a+b}$. Let $u = a+b$. Applying the comultiplication gives

$$\sum_{u_1+u_2=u} \binom{u}{a} Z_{u_1} \otimes Z_{u_2} \tag{3.17}$$

On the other hand

$$\sum_{\substack{a_1+a_2=a \\ b_1+b_2=b}} \text{coeff}(s^{a_1}t^{b_1}) \; \text{coeff}(s^{a_2}t^{b_2}) = \sum_{\substack{a_1+a_2=a \\ b_1+b_2=b}} \binom{a_1+b_1}{a_1} Z_{a_1+b_1} \binom{a_2+b_2}{a_2} Z_{a_2+b_2} \quad (3.18)$$

Take any $u_1$, $u_2 = u - u_1$. Then the coefficient of $Z_{u_1} Z_{u_2}$ in (3.18) is equal to

$$\sum_{a_1} \binom{u_1}{a_1} \binom{u_2}{a-a_1} = \binom{u_1+u_2}{a} = \binom{u}{a} \quad (3.19)$$

Where the binomial coefficient identity (3.19) follows from looking at the coefficient of $t^a$ in

$$(1+t)^{u_1}(1+t)^{u_2} = (1+t)^{u_1+u_2}$$

which proves the lemma. From the dual point of view things (when applicable, which is certainly the case here), things are much easier. If $d(t)$ is a curve in $H$ and $H \to \mathbf{Z}[[t]]$ is the corresponding morphisms of algebras, then $d(t+s)$ is the 2-curve correponding to the composed morphism

$$H \to \mathbf{Z}[[t]] \xrightarrow{t \mapsto s+t} \mathbf{Z}[[s,t]].$$

The next bit of preparation concerns the noncommutative Newton primitives

$$P_n(Z) = \sum_{r_1+\cdots+r_k=n} (-1)^{k+1} r_k Z_{r_1} Z_{r_2} \cdots Z_{r_k}, \quad r_i \in \mathbf{N} = \{1,2,\cdots\} \quad (3.20)$$

(Note that these differ by a sign factor $(-1)^{n+1}$ from the slightly more often used Newton primitives, $Z_1$, $Z_1^2 - 2Z_2$, $Z_1^3 - 2Z_1Z_2 - Z_2Z_1 + 3Z_3, \cdots$.). The Newton primitives (3.20) satisfy the recursion relation

$$P_n(Z) = nZ_n - Z_{n-1}P_1(Z) - Z_{n-2}P_2(Z) - \cdots - Z_1 P_{n-1}(Z) \quad (3.21)$$

Now, over *NSymm*, consider the 2-curve $Z(s)^{-1}Z(t)^{-1}Z(s+t)$.

3.22. **Theorem** (Second isobaric decomposition theorem). There are unique homogeneous noncommutative polynomials $N_{u,v}(Z) \in NSymm$ such that

$$Z(s)^{-1}Z(t)^{-1}Z(s+t) = \prod_{\substack{a,b \in \mathbf{N} \\ \gcd(a,b)=1}} (1 + N_{a,b}(Z)s^a t^b + N_{2a,2b}(Z)s^{2a}t^{2b} + \cdots). \quad (3.23)$$

Moreover

(i) $N_{u,v}(Z) = \binom{u+v}{u} Z_{u+v} + \text{(terms of length} \geq 2)$ (3.23)

(ii) $N_{u,v}(Z)$ is homogeneous of weight $u+v$ (3.24)

(iii) For each $a,b \in \mathbf{N}^2$, $\gcd(a,b) = 1$,

$$1 + N_{a,b}(Z)s^a t^b + N_{2a,2b}(Z)s^{2a}t^{2b} + \cdots \quad (3.25)$$

is a 2-curve.

(iv) For each $n \geq 2$, $N_{1,n-1}(Z) = P_n(Z)$ (3.26)

Proof. The situation is very like the one of the first decomposition theorem above. Again putting $s$ or $t$ equal to zero in the left hand side of (3.23) gives 1, so that there are no pure powers of $s$ and $t$ in (3.23), and, again, for every pair $(u,v) \in \mathbf{N}^2$ there is precisely one of the factors on the right hand side of (3.23) in which $s^u t^v$ occurs. Proceeding as, before, i.e. bring $Z(s)$ and $Z(t)$ over to the right hand side, and compare coefficients, one finds a recursion formula for the $N_{u,v}(Z)$

$$N_{u,v}(Z) = \binom{u+v}{u} Z_{u+v} - \sum_{\substack{u_0 + \cdots + u_k = u \\ v_0 + \cdots + v_k = v \\ u_i, v_i \geq 1, \text{ for } i \geq 1 \\ u_0 + v_0 > 0 \text{ or } k \geq 2}} Z_{v_0} Z_{u_0} N_{u_1, v_1}(Z) \cdots N_{u_k, v_k}(Z) \quad (3.27)$$

where the product is again an ordered one; i.e. (3.13) must hold. Using (3.27) as a definition it follows that (3.23) holds. This takes care of existence and uniqueness. Properties (3.23) and (3.24) follow immediately from the recursion formula (3.27). Note that there are at least two factors in each of the terms under the sum sign in (3.27) so that all terms there have weight less than $u + v$.

The proof that the (3.25) are curves is exactly the same as in the case of the first isobaric decomposition theorem, using that the left hand side of (3.23) is a 2-curve; see Lemma 3.17.

Finally, take $u = 1 = v$ in (3.27). The only possible term in this case under the sum sign on the right has $k = 0$, $u_0 = 1$, $v_0 = 1$ and thus

$$N_{1,1}(Z) = 2Z_2 - Z_1^2$$

Now let $v \geq 2$. Then the recursion formula (3.27) gives

$$N_{1,v} = \binom{v+1}{1} Z_{v+1} - Z_1 N_{1,v-1} - Z_2 N_{1,v-2} - \cdots - Z_{v-1} N_{1,1} - Z_v Z_1 \quad (3.28)$$

Noting that $P_1(Z) = Z_1$ it follows with induction that this is the same recursion formula as for the $P_n(Z)$ (with $n = v - 1$, see (3.21) above), proving the last statement of the theorem.

3.29. Remarks. Using a different ordering in the ordered products occurring in (3.7) and (3.23), and using slightly different 'commutation formulae' one obtains different versions of the $L_{u,v}(X,Y)$ and $N_{u,v}(Z)$ with nice symmetry properties. This does not matter for the purposes of the present paper but probably deserves further exploration for future applications and calculations.

The first thing one needs is an ordering on the set $J = \{(a,b) \in \mathbf{N} \times \mathbf{N}: \gcd(a,b) = 1\}$ with the property that $(a,b) >_{swl} (a',b') \Leftrightarrow (b',a') >_{swl} (b,a)$. There are many such. Here is one that fits the present circumstances rather well. Divide $J$ into three parts as follows:

$$J = J_- \cup \{(1,1)\} \cup J_+ \quad (3.30)$$

where $J_- = \{(a,b) \in J: a > b\}$, $J_+ = \{(a,b) \in J: a < b\}$ Now on $J_-$ take the ordering 'weight first and lexicographic afterwards, and on $J_+$ take the reverse ordering; i.e $(a,b) >_{swl} (c,d)$ for $(a,b),(c,d) \in J_+$ if and only if $(d,c) >_{swl} (b,a)$ in $J_-$. Further set $J_- <_{swl} (1,1) <_{swl} J_+$ Thus restricting to weight $\leq 7$, the resulting ordering is

$$(2,1) <_{swl} (3,1) <_{swl} (3,2) <_{swl} (4,1) <_{swl} (5,1) <_{swl} (4,3) <_{swl} (5,2) <_{swl} (6,1) <_{swl} \cdots <_{swl} (1,1)$$
$$\cdots <_{swl} (1,6) <_{swl} (2,5) <_{swl} (3,4) <_{swl} (1,5) <_{swl} (1,4) <_{swl} (2,3) <_{swl} (1,3) <_{swl} (1,2)$$

The suffix 'wsl' is supposed to be an acronym for 'symmetric weight first lexicographic after'. There is a first element, viz (2,1) and a last element, viz. (1,2), and there are inifitely many elements between the first one and the 'middle one', (1.1), and between the middle one and the last one. In this it is a somewhat unusual order . In each equal weight segment the order is lexicographic.

Now there are a bi-isobaric decompositions

$$X(s)^{-1}Y(t)X(s)Y(t)^{-1} = \sum_{\gcd(a,b)=1} (1 + L_{a,b}(X,Y)s^a t^b + L_{2a,2b}(X,Y)s^{2a}t^{2b} + \cdots) \qquad (3.31)$$

$$Z(s)^{-1}Z(s+t)Z(t)^{-1} = \sum_{\substack{a,b \in \mathbf{N} \\ \gcd(a,b)=1}} (1 + N_{a,b}(Z)s^a t^b + N_{2a,2b}(Z)s^{2a}t^{2b} + \cdots) \qquad (3.32)$$

where now the ordering is the one just defined, i.e. the $<_{swl}$-ordering, with corresponding recursion forrmulas

$$L_{u,v}(X,Y) = -[X_u, Y_v] - \sum_{\substack{u_0,v_0 \geq 0; k \geq 1 \\ u_i,v_i \geq 1, i=1,\cdots,k \\ u_0+u_1+\cdots+u_k=u \\ v_0+v_1+\cdots+v_k=v \\ u_0+v_0 \geq 1 \text{ or } k \geq 2}} X_{u_0} L_{u_i,v_i}(X,Y) Y_{v_0} \qquad (3.33)$$

$$N_{u,v}(Z) = \binom{u+v}{u} Z_{u+v} - \sum_{\substack{u_0+\cdots+u_k=u \\ v_0+\cdots+v_k=v \\ u_i,v_i \geq 1, \text{ for } i \geq 1 \\ u_0+v_0>0 \text{ or } k \geq 2}} Z_{u_0} N_{u_1,v_1}(Z) \cdots N_{u_k,v_k}(Z) Z_{v_0} \qquad (3.34)$$

Let $\rho$ on $\mathbf{Z}\langle X;Y \rangle$ or $\mathbf{Z}\langle Z \rangle$ be the anti-isomorphism of algebras that reverses the order of multiplication. Thus, e.g. $\rho(Z_1 Z_3 Z_2) = Z_2 Z_3 Z_1$, $\rho([X_i, Y_j]) = [Y_j, X_i] = -[X_i, Y_j]$. Then there are the symmetry properties

$$L_{v,u}(Y,X) = \rho(L_{u,v}(X,Y)) \qquad (3.35)$$

$$N_{v,u}(Z) = \rho(N_{u,v}(Z)) \qquad (3.36)$$

In this case

$$N_{k,1} = P_{k+1}(Z), \quad N_{1,k}(Z) = Q_{k+1}(Z)$$

where the $Q_n(Z)$ are the second family of 'power sum primitives' defined by the recursion formula

$$Q_n(Z) = nZ_n - Q_1(Z)Z_{n-1} - Q_2(Z)Z_{n-2} - \cdots - Q_{n-1}(Z)Z_1 \qquad (3.37)$$

One reason that the ordering $<_{swl}$ is esthetically nice is that one can insert the element (1,0) at the beginning, before (2,1), and (0,1) at the end, after (1,2). Now write

$$L_{u,0}(X,Y) = X_u, \quad L_{0,v}(X,Y) = Y_v, \quad N_{u,0}(Z) = Z_u = N_{0,u}(Z)$$

and then the products in the sums in (3.33) and (3.34) are all according to the ordering

$$(1,0) <_{swl} (2,1) <_{swl} (3,1) <_{swl} (3,2) <_{swl} (4,1) <_{swl} (5,1) <_{swl} (4,3) <_{swl} \cdots <_{swl} (1,1)$$
$$\cdots <_{swl} (3,4) <_{swl} (1,5) <_{swl} (1,4) <_{swl} (2,3) <_{swl} (1,3) <_{swl} (1,2) <_{swl} (0,1)$$

There is also a second symmetry property, obtained by taking inverses of (3.31) and (3.32) involving the antipode. Here also an order like $<_{wsl}$ is important because taking inverses reverses the product order and then switching $s$ and $t$ restores it again.

### 4. The primitives of *NSymm*.

The construction of the basis $P_\alpha$, $\alpha \in LYN$ was already given in the introduction. It is repeated here for convenience of reference and with more detail, proofs, and explanations.

At this stage we have the naturally given (supreme) curve

$$z, \quad z(i) = Z_i; \quad z(t) = 1 + Z_1 t + Z_2 t^2 + \cdots \qquad (4.1)$$

in *NSymm* (using all three notations so far used); and, using theorem 3.5 (The Ditters-Shay bi-isobaric decomposition theorem) and Lemma 2.2.4, a large number of curves

$$c_{a,b}, \quad c_{a,b}(i) = L_{ia,ib}(X,Y), \quad c_{a,b}(t) = 1 + L_{a,b}(X,Y)t + L_{2a,ab}(X,Y)t^2 + \cdots \qquad (4.2)$$

in $2NSymm$, one for each pair of positive integers $(a,b), a,b \in \mathbf{N}^2, \mathbf{N} = \{1,2,3,\cdots\}$ with $\gcd(a,b) = 1$. Substituting curves in curves in a suitable way gives (potentially) new curves. Below there is a systematic procedure of doing this that leads to a basis of $\mathrm{Prim}(NSymm)$.

Let $\alpha$ be a composition (word) over the positive integers

$$\alpha = [a_1, a_2, \cdots, a_m], \quad a_i \in \mathbf{N} \qquad (4.3)$$

Its weight and length are respectively

$$\mathrm{wt}(\alpha) = a_1 + \cdots + a_m, \quad \mathrm{lg}(\alpha) = m \qquad (4.4)$$

The corresponding noncommutative monomials

$$Z_\alpha = Z_{a_1} Z_{a_2} \cdots Z_{a_m} \qquad (4.5)$$

have the same weight and length. The empty composition and the monomial 1 have length and weight zero.

For use in a minute or so, here is the wll-ordering on compositions and monomials:

$$\alpha <_{wll} \beta \quad \begin{array}{l} wt(\alpha) < wt(\beta) \\ \text{or } (wt(\alpha) = wt(\beta) \text{ and } \mathrm{lg}(\alpha) < \mathrm{lg}(\beta)) \\ \text{or } (wt(\alpha) = wt(\beta) \text{ and } \mathrm{lg}(\alpha) = \mathrm{lg}(\beta) \text{ and } \alpha <_{lexico} \beta) \end{array} \qquad (4.6)$$

The index 'wll' is supposed to be a mnemonic for 'weight first, than length, and thereafter lexocographic'.

A (nonempty) composition $\alpha$ is Lyndon it it is lexicographically smaller than all its proper tails, $[a_i, \cdots, a_m]$, $i = 2, \cdots m$. The set of Lyndon words is denoted $LYN$. To each Lyndon word there are associated three things

- a positive integer $g(\alpha) = \gcd(a_1, a_2, \cdots a_m)$ (4.7)
- a curve $d_\alpha$
- a primitive $P_\alpha = P_{g(\alpha)}(d_\alpha)$ (4.8)

where the $P_n(Z)$ are the Newton primitives as defined by (3.20). That $P_\alpha$ is indeed a primitive follows from the fact that the $P_n(Z)$ are primitives (see e.g. [3]) and that, for any curve $d$, $P_n(d)$ is the image of $P_n(Z)$ under the Hopf algebra morphism corresponding to the curve $d$, and hence a primitive.

It remains to give the recursive construction of the $d_\alpha$. This goes by induction on the length of $\alpha$.

- For $\lg(\alpha) = 1$, i.e. $\alpha = [n]$ for some $n$, $d_{[n]} = z$
- For $\lg(\alpha) > 1$, let $\alpha'' = [a_i, \cdots, a_m]$ be the lexicographically smallest proper tail of $\alpha$, and let $\alpha' = [a_1, \cdots, a_{i-1}]$. Then $\alpha', \alpha''$ are both Lyndon. (Indeed, $\alpha''$ is Lyndon because it is the lexicographically *smallest* proper tail of $\alpha$; further, if $[a_j, \cdots, a_{i-1}]$ were lexicographicaly larger than $\alpha'$, then $[a_j, \cdots, a_m]$ would be lexicographically larger than $[a_1, \cdots, a_m] = \alpha$, which is not the case; so $\alpha'$ is also Lyndon). The concatenation factorization $\alpha = \alpha' \alpha''$ is called the canonical factorization of the Lyndon word $\alpha$. The curve $d_\alpha$ is now defined by

$$d_\alpha = c_{g(\alpha')/g(\alpha), g(\alpha'')/g(\alpha)}(d_{\alpha'}, d_{\alpha''}) \qquad (4.9)$$

Note that $d_{r\alpha} = d_\alpha$ for any $r \in \mathbf{N}$, where $r\alpha = [ra_1, ra_2, \cdots, ra_m]$.

4.10. **Proposition.** The wll-smallest term in $d_\alpha(i)$ is $Z_{ig(\alpha)^{-1}\alpha}$ and it occurs with coefficient 1.

Proof. This is proved by induction on length. The case of length 1 is obvious as $d_{[n]} = z$ for all $n$. Now let $\alpha = \alpha' \alpha''$ be the canonical factorization of $\alpha$. First observe that by theorem 3.5 (ii) (formula (3.9)), there are no monomials in the $L_{u,v}(X,Y)$ that consist purely of $X$'s or purely of $Y$'s. Using theorem (3.5) (i) and calculating modula $\lg(\alpha) + 1$ is follows that

$$d_\alpha(i) = L_{ig(\alpha)^{-1}g(\alpha'), ig(\alpha)^{-1}g(\alpha'')}(d_{\alpha'}, d_{\alpha''}) \equiv [d_{\alpha'}(ig(\alpha)^{-1}g(\alpha')), d_{\alpha''}(ig(\alpha)^{-1}g(\alpha''))] \qquad (4.11)$$

By induction the wll-smallest terms of the two expressions in the commutator on the right of (4.11) are, respectively

$$Z_{ig(\alpha)^{-1}\alpha'}, \quad Z_{ig(\alpha)^{-1}\alpha''}$$

and thus the wll-smallest term in $d_\alpha(i)$ is either

$$Z_{ig(\alpha)^{-1}\alpha'} Z_{ig(\alpha)^{-1}\alpha''} = Z_{ig(\alpha)^{-1}\alpha}$$

with coefficient 1, or $Z_{ig(\alpha)^{-1}\alpha''} Z_{ig(\alpha)^{-1}\alpha'}$ with coefficient -1. However, $ig(\alpha)^{-1}\alpha'' \cdot ig(\alpha)^{-1}\alpha'$ is lexicographically larger than $ig(\alpha)^{-1}\alpha' \cdot ig(\alpha)^{-1}\alpha''$ because $\alpha$ is Lyndon. Here '·' denotes concatenation.

4.12. **Corollary**. The wll-smallest term in the primitive $P_\alpha$ is $g(\alpha)Z_\alpha$ for all $\alpha \in LYN$.

Proof. this follows from the defining formula (4.8) and the observation that $P_n(Z) \equiv nZ_n$ modulo (length $\geq 2$).

4.13. **Corollary**. The $P_\alpha$, $\alpha \in LYN$, form a basis of Prim(*NSymm*) over the rationals.

Proof. Let $\beta_n$ be the number of Lyndon words of weight $n$. The $P_\alpha$ with $\alpha$ Lyndon of weight $n$ are homogeneous of weight $n$ and they are linearly independent because when tested against the monomials $Z_\alpha$, $\alpha \in LYN_n = \{\alpha \in LYN: \text{wt}(\alpha) = n\}$ they form a triangular matrix (with $g(\alpha)$'s on the diagonal (ordering $LYN_n$ by the wll-ordering)). Now over the rationals *NSymm* is isomorphic to

$$\mathcal{U} = \mathbf{Z}\langle U_1, U_2, \cdots \rangle, \quad \mu(U_n) = 1 \otimes U_n + U_n \otimes 1$$

(see (1.6)); and the isomorphism is degree preserving. The Lie algebra of primitives of $\mathcal{U}$, Prim($\mathcal{U}$) is the free Lie algebra generated by the $U_1, U_2, \cdots$ and has a homogeneous basis indexed by Lyndon words, $Q_\alpha$, $\alpha \in LYN$, $\text{wt}(Q_\alpha) = \text{wt}(\alpha)$. In particular $\dim(\text{Prim}(\mathcal{U})_n) = \beta_n$. It follows that the isomorphism maps the space spanned by the $P_\alpha$, $\alpha \in LYN$ onto Prim($\mathcal{U}$)$_n$, proving the corollary.

The collection of independent primitives $P_\alpha$, $\alpha \in LYN$ is very well behaved with respect to the Verschiebung Hopf algebra morphisms $\mathbf{V}_r$.

4.14. **Theorem**. For all $\alpha \in LYN$

$$\mathbf{V}_r(P_\alpha) = \begin{cases} rP_{r^{-1}\alpha} & \text{if } r \text{ divides } g(\alpha) \\ 0 & \text{if } r \text{ does not divide } g(\alpha) \end{cases} \tag{4.15}$$

Proof. First consider the case of length 1, $P_{[n]} = P_n(Z)$, for which there is an explicit formula, viz

$$P_n(Z) = \sum_{r_1 + \cdots + r_k = n} (-1)^{k+1} r_k Z_{r_1} Z_{r_2} \cdots Z_{r_k}, \quad r_i \in \mathbf{N} = \{1,2,\cdots\} \tag{4.16}$$

Now $\mathbf{V}_r$ applied to one of the monomials making up (4.16) gives zero unless each of the $r_i$, $i = 1, \cdots, k$ is divisible by $r$. This can happen only if $r$ divides $n = r_1 + \cdots + r_k$, proving that

$$\mathbf{V}_r P_n(Z) = 0 \text{ if } r \text{ does not divide } n \tag{4.17}$$

If $r$ does divide all the $r_i$ then the result of applying $\mathbf{V}_r$ to $(-1)^{k+1} r_k Z_{r_1} Z_{r_2} \cdots Z_{r_k}$ is

$$r(-1)^{k+1} u_k Z_{u_1} Z_{u_2} \cdots Z_{u_k}, \quad ru_j = r_j, \quad u_1 + \cdots + u_k = r^{-1}n$$

and summing this over all possibilities precisely gives $rP_{r^{-1}n}(Z)$. Thus

$$\mathbf{V}_r P_n(Z) = rP_{r^{-1}n}(Z) \text{ if } r \text{ divides } n \tag{4.18}$$

Writing (4.17) and (4.18) out a bit more gives also

$$P_n(\underbrace{0,\cdots,0}_{r-1},Z_1,\underbrace{0,\cdots,0}_{r-1},Z_2,\underbrace{0,\cdots,0}_{r-1},Z_3,\cdots) = \begin{matrix} 0 \text{ if } r \text{ does not divide } n \\ rP_{r^{-1}n}(Z_1,Z_2,Z_3,\cdots) \text{ if } r \text{ divides } n \end{matrix} \quad (4.19)$$

Now let $\alpha$ be any Lyndon word. Then $P_\alpha = P_{g(\alpha)}(d_\alpha(1),\cdots,d_\alpha(g(\alpha)))$. Now by theorem 3.5 (iii), the $c_{a,b}$ are **V**-curves, and the recursive construction of $d_\alpha$ starts with the **V**-curve $z$ so $d_\alpha$ is a **V**-curve by proposition 2.4.3. Thus, using (4.19),

$$\mathbf{V}_r P_{g(\alpha)}(d_\alpha(1),d_\alpha(2),\cdots) = P_{g(\alpha)}(\mathbf{V}_r d_\alpha(1),\mathbf{V}_r d_\alpha(2),\cdots) =$$

$$= P_{g(\alpha)}(\underbrace{0,\cdots,0}_{r-1},d_\alpha(1),\underbrace{0,\cdots,0}_{r-1},d_\alpha(2),\cdots) = \begin{matrix} 0 \text{ if } r \text{ does not divide } n \\ rP_{r^{-1}g(\alpha)}(d_\alpha(1),d_\alpha(2),\cdots) = rP_{r^{-1}\alpha} \text{ if } r \text{ divides } n \end{matrix}$$

proving the theorem.

In the proof below of the main theorem, the theorem that the $P_\alpha$, $\alpha \in LYN$, form a basis over the integers for Prim(*NSymm*), the following lemma is used.

4.17. Lemma. Let $p$ be a prime number, then

$$\mathbf{V}_p \text{Prim}(NSymm) \subset p\text{Prim}(NSymm)$$

The proof of this involves the graded dual Hopf algebra of *NSymm*, the Hopf algebra *QSymm* of quasisymmetric functions. This Hopf algebra will again turn up later in this paper. For a good deal of information on the Hopf algebra of symmetric functions see [6] and the references therein. Here is a brief description of only some of its structure. Being the graded dual of *NSymm*, which has the monomials $Z_\alpha$, $\alpha$ running over all words over the positive integers, as an Abelian group basis one can take as a basis for *QSymm* those words themselves with the duality pairing

$$\langle,\rangle : NSymm \times QSymm \to \mathbf{Z}, \quad \langle Z_\alpha, \beta \rangle = \delta_\alpha^\beta \quad (4.18)$$

(Kronecker delta). The multiplication on *QSymm* induced by the pairing (4.18) is the overlapping shuffle multiplication. It can be described as follows. Let $\alpha = [a_1,\cdots,a_m]$ and $\beta = [b_1,\cdots b_n]$ be two compositions. For each $k$, $0 \le k \le \min(m,n)$ take a word of length $m + n - k$ with sofar empty slots. Divide the slots in three disjoint sets $A$, $B$, $O$ of sizes $m-k$, $n-k$, $k$ respectively. There are $\binom{m+n-k}{m-k,n-k,k}$ ways of doing this. In the slots of $A \cup O$ put the $a$'s in their original order; in the slots of $B \cup O$ put the $b$'s in their original order. In the slots of $O$ there is both an $a$ and a $b$ entry; add those. The product $\alpha\beta = \alpha \times_{osh} \beta$ is the sum of all the words thus obtained. Thus the product has $\sum_k \binom{m+n-k}{m-k,n-k,k}$ terms. For instance

$$[a][b_1,b_2] = [a,b_1,b_2] + [b_1,a,b_2] + [b_1,b_2,a] + [a+b_1,b_2] + [b_1,a+b_2]$$

The comultiplication on *QSymm* is 'cut'

$$\mu([a_1,\cdots,a_m]) = [] \otimes [a_1,\cdots,a_m] + \sum_{i=1}^{m-1}[a_1,\cdots,a_i] \otimes [a_{i+1},\cdots,a_m] + [a_1,\cdots,a_m] \otimes []$$

where [] is the empty word which is the unit element of *QSymm*.

It is useful to have a description in one go of the overlapping shuffle product of possibly more than two factors. For a matrix $M$ with entries from $\mathbf{N} \cup \{0\}$ let $\alpha_c(M)$ denote the compostion of the column sums of the entries of $M$. Let $\alpha_1, \alpha_2, \cdots, \alpha_k$ be $k$ compositions. Now, consider all matrices $M$ of $k$ rows that are as follows:
- the $i$-th row of $M$ consists of zeros and the entries of $\alpha_i$ in their original order,
- there are no columns consisting entirely of zeros.

Then the product $\alpha_1 \alpha_2 \cdots \alpha_k$ is equal to

$$\alpha_1 \alpha_2 \cdots \alpha_k = \sum_M a_c(M) \tag{4.19}$$

where $M$ runs over all different matrices satisfying the two conditions just described.

The endomorphisms of Hopf algebras of *QSymm* dual to the Verschiebung endomorphisms of *NSymm* are the Frobenius morphisms

$$\mathbf{f}_n([a_1, a_2, \cdots, a_m]) = [na_1, na_2, \cdots, na_m] \tag{4.20}$$

This is immediate from the duality pairing (4.18).

4.21. **Lemma.** For each prime number $p$

$$\mathbf{f}_p \alpha \equiv \alpha^p \mod p \tag{4.22}$$

Proof. Let $\alpha = [a_1, \cdots, a_m]$. Consider a matrix $M$ such as in the description of the multiple overlapping shuffle product of compositions given above and consider a (nonidentity) cyclic permutation of its rows. Unless all the rows are identical the result is a different matrix (because $p$ is a prime number and the cyclic group of of $p$ elements has no proper subgroups and the stabilizer of an element of a set acted on by a group must be a subgroup of that group). The column sums of a matrix and a row permuted version of that matrix are the same. There is only one matrix with all rows identical and this one gives the term $\mathbf{f}_p \alpha = [pa_1, \cdots, pa_m]$. For all other matrices of the required properties the $p$ cyclic permutations are all different and thus give rise to terms of the form $p(\text{something})$. This proves the lemma.

4.23. **Proof of Lemma 4.17.** Let $P$ be a primitive element of *NSymm* and let $\alpha$ be any composition. Then (modulo $p$)

$$\langle V_p P, \alpha \rangle = \langle P, \mathbf{f}_p \alpha \rangle = \langle P, \alpha^p \rangle = \langle \mu_p(P), \underbrace{\alpha \otimes \cdots \otimes \alpha}_{p} \rangle =$$

$$= \langle \underbrace{1 \otimes \cdots \otimes 1}_{p-1} \otimes P + \underbrace{1 \otimes \cdots \otimes 1}_{p-2} \otimes P \otimes 1 + \cdots + P \otimes \underbrace{1 \otimes \cdots \otimes 1}_{p-1}, \underbrace{\alpha \otimes \cdots \otimes \alpha}_{p} \rangle = 0$$

where $\mu_p$ is the $p$-fold coproduct ($\mu_2 = \mu$, $\mu_3 = (id \otimes \mu)\mu$, $\mu_4 = (id \otimes id \otimes \mu)\mu_3 \cdots$).

Things are now ready for the main theorem and its proof.

4.24. **Theorem.** The $P_\alpha$, $\alpha \in LYN$ form a basis (over the integers) for Prim(*NSymm*).

Proof. Prim(*NSymm*) is a graded subgroup of *NSymm*; let Prim(*NSymm*)$_n$ be the piece of weight $n$. By Corollary 4.13, the $P_\alpha$ with $\alpha$ Lyndon of weight $n$, i.e. $\alpha \in LYN_n$, span a full rank subgroup of Prim(*NSymm*)$_n \subset$ *NSymm*$_n$. It therefore only remains to show that that it is a pure subgroup of *NSymm*$_n$.

( A pure subgroup of a free Abelian group $A$ of finite rank is a subgroup $G$ such that if for a prime number $p$ and $a \in A$, $pa \in G$, then $a \in G$. Now $\mathrm{Prim}(NSymm)_n \subset NSymm_n$ is a pure subgroup, because *NSymm* is torsion free. Finally if $G_1 \subset G_2$, are subgroups of the torsion free Abelian finite rank group $A$, $G_1$ is pure, and $G_1$ and $G_2$ have the same rank than they are equal).

The proof that (Subgroup generated by the $P_\alpha$, $\alpha \in LYN_n$) $\subset NSymm_n$ is pure), is done by induction on weight. The case of weight 1 is trivial

$$NSymm_1 = \mathbf{Z}Z_1 = \mathrm{Prim}(NSymm)_1 = \mathbf{Z}P_{[1]}$$

(and the case of weight two is easy with $\mathrm{Prim}(NSymm)_2$ of rank 1 and $\mathbf{Z}P_{[2]} = \mathbf{Z}(2Z_2 - Z_1^2)$ manifestly a pure subgroup of $NSymm_2 = \mathbf{Z}Z_2 \oplus \mathbf{Z}Z_1^2$; one can also easily deal with the cases weight 3,4,5 by direct inspection; after that things become rapidly too complicated for direct hand calculation).

So take an $n \geq 2$ and suppose that pureness has been proved for all smaller weights. When tested against the Lyndon monomials of weight $n$, the $P_\alpha$ with $\alpha$ of weight $n$, yield a triangular matrix with the $g(\alpha)$'s on the diagonal. It follows that the only primes to worry about are the divisors of the determinant $\prod_{\alpha \in LYN_n} g(\alpha)$. These are all divisors of $n$. So let

$$Q = \sum_{\alpha \in LYN_n} x_\alpha P_\alpha \in pNSymm_n \quad (4.25)$$

where the prime number may be assumed to be a divisor of $n$. There is to prove that all the integer coefficients $x_\alpha$ are divisible by $p$. By the nature of the primitiveness condition $p^{-1}Q = Q \in \mathrm{Prim}(NSymm)_n$ (as *NSymm* is torsion free). Hence, by Lemma 4.17, $\mathbf{V}_p Q \in p^2 NSymm$. Now apply $\mathbf{V}_p$ to equation (4.25) to get (using Theorem 4.14)

$$\mathbf{V}_p(\sum_{\alpha \in LYN_n} x_\alpha P_\alpha) = \sum_{\substack{\alpha \in LYN_n \\ p \text{ divides } g(\alpha)}} p x_\alpha P_{p^{-1}\alpha} = \mathbf{V}_p Q \in p^2 NSymm_n$$

So that

$$\sum_{\substack{\alpha \in LYN_n \\ p \text{ divides } g(\alpha)}} x_\alpha P_{p^{-1}\alpha} \in pNSymm_n \quad (4.26)$$

Now the $p^{-1}\alpha$ with $\alpha \in LYN_n$ and $p | g(\alpha)$ are precisely the elements of $LYN_{p^{-1}n}$. And thus, by induction, we see from (4.26) that the coefficients $x_\alpha$ in (4.25) for which $p$ divides $g(\alpha)$ are divisible by $p$. Thus these terms can be moved to the aother side of (4.25) to give a relation

$$\sum_{\substack{\alpha \in LYN_n \\ \gcd(p, g(\alpha))=1}} x_\alpha P_\alpha \in pNSymm_n \quad (4.27)$$

Let $\langle P_\alpha, \beta \rangle$ be the coefficient of the monomial $Z_\beta$ in $P_\alpha$. (This is just the duality pairing mapping (4.18)). Then the matrix ( $\langle P_\alpha, \beta \rangle$ ) for $\alpha, \beta \in LYN_n$ and $\gcd(p, g(\alpha)) = 1 = \gcd(p, g(\beta))$ is triangular with determinant prime to $p$ (Corollary 4.12). Together with (4.27) this proves that the prime number $p$ divides all the coefficients $x_\alpha$ in (4.27). This finishes the proof.

4.28. Comment. In the proof above rather heavy use has been made of the notion of a **V**-curve (via Theorem 4.14). This can be avoided, as in [12], at the price of an additional double induction (with respect to the higher commutator ideals of *NSymm* and with respect to the wll-ordering). The essential thing being to show that the top term (wll-ordering) in (4.25) does not go to zero under $V_p$ (if $p$ divides its index).

4.29. Corollary (of Theorem 4.24 together with Theorem 4.14).

$$\mathbf{V}_r \text{Prim}(NSymm) = r\text{Prim}(NSymm) \tag{4.30}$$

This has been proved in [1] under the assumption that the graded dual of *NSymm* is free over the integers.

The $P_\alpha$, $\alpha \in LYN$ for which $g(\alpha) = 1$ are by construction the first term of a (curve) divided power sequence, viz $d_\alpha$. For the $P_\alpha$, $\alpha \in LYN$ with $g(\alpha) > 1$, there is also a (natural) curve of which it is the first term. Indeed, consider the curve

$$1 + N_{1,g(\alpha)-1}(d_\alpha)t + N_{2,2g(\alpha)-2}(d_\alpha)t^2 + N_{3,3g(\alpha)-3+}(d_\alpha)t^3 + \cdots \tag{4.31}$$

where the $N_{u,v}(Z)$ are the noncommutative polynomials in *NSymm* defined by the second isobaric decomposition theorem 3.22. Thus

4.32. Theorem. The Lie algebra of primitives of *NSymm* has a basis consisting of homogenous elements that each extend to an isobaric curve. There are $\beta_n$ basis elements of weight *n*.

Here an isobaric curve over a homogeneous primitive of weight $n$ is a curve $d = (1, P, d(2), d(3), \cdots)$ such that $\text{wt}(d(k)) = kn$, and $\beta_n = \#LYN_n$.

Theorem 4.32 has a most important consequence

4.33. Theorem. The graded dual Hopf algebra of *NSymm*, the Hopf algebra *QSymm* of quasisymmetric functions, is free (as a commutative algebra) over the integers with $\beta_n$ generators of weight $n$.

There is a proof of the implication 4.32 ⇒ 4.33 in [12]. One can also adapt the proof in [13], Chapter XIII, p.273ff of the same result for Hopf algebras over a field of characteristic zero. For completeness sake there is a simplified proof (which has ingredients from both) in the appendix, adapted to this special case.

For another, completely different, proof of the freeness theorem 4.33, see [5].

**Appendix**

A.1. *The graded cocommutative cofree coalgebra over a graded set*.
Let $S$ be a graded set; i.e. each element of $S$ has a weight attached to it that is a positive integer. Only those graded sets are considered for which there are only finitely many elements of each weight. Let $\gamma_n$ be the number of elements of $S$ of weight $n$.

The graded commutative free algebra over **Z** on $S$ is the algebra of polynomials

$$\mathrm{CF}(S) = \mathbf{Z}[x_s : s \in S], \quad \mathrm{wt}(x_s) = \mathrm{wt}(s) \tag{A.1.1}$$

Give the set $S$ a total ordering; it is natural to do this in such a way that the elements of weight 1 come first, then those of weight 2, etc.

If $f : S \to \mathbf{N} \cup \{0\}$ is a function of finite support to the nonnegative integers, let

$$x_f = \prod_{s \in S, f(s) > 0} x_s^{f(s)} \tag{A.1.2}$$

where the product is ordered by the ordering of $S$. These monomials form an Abelian group basis for $\mathrm{CF}(S)$. The number of monomials in the $x_s$ of weight $n$ is the coefficient of $t^n$ in the power series expression

$$\frac{1}{(1-t)^{\gamma_1}(1-t^2)^{\gamma_2}(1-t^3)^{\gamma_3}\cdots} \tag{A.1.3}$$

The graded cofree cocommutative coalgebra over $S$, CCoF($S$), is the graded dual of CF($S$). The rank of its homogeneous part of weight $n$ is of course again the coefficient of $t^n$ in (A.1.3). By duality it is easy to figure out the comultiplication on CCoF($S$). By definition a basis (as an Abelian group) is given by the finite support functions $g: S \to \mathbf{N} \cup \{0\}$ and the duality pairing is

$$\langle g, x_f \rangle = \delta^g_f \qquad (A.1.4)$$

(Kronecker delta). It follows immediately that the comultiplication is

$$\mu(f) = \sum_{f_1 + f_2 = f} f_1 \otimes f_2 \qquad (A.1.5)$$

The counit takes each $f$ that is not identically zero to zero and the identically zero function on $S$ to 1. The graded cofree cocommutative coalgebra over $S$, CCoF($S$), has an appropriate universality property as follows. Let $M(S)$ be the graded free module spanned by $S$ and $M(S)^*$ its graded dual module. Then CCoF($S$) comes with a natural morphism of Abelian groups CCoF($S$) $\to M(S)^*$ (the dual of the natural inclusion $M(S) \to CF(S)$), and satisfies the following universality property: for each graded cocommutative coalgebra $C$ and morphism of graded Abelian groups $C^+ \to M(S)^*$ there is a unique lift $C \to$ CCoF($S$) that is a morphism of coalgebras. Here $C^+ = Ker(\varepsilon)$ where $\varepsilon$ is the counit of CCoF($S$). This universality property will not be used here.

For a construction via the tensor power sum $\mathbf{Z} \oplus M \oplus M^{\otimes 2} \oplus M^{\otimes 3} \oplus \cdots$, where $M$ is short for $M(S)$, see [13], Chapter XII, p.243ff.

A.2. *Proof of the freeness theorem* (see Theorem 4.33 above).
For each Lyndon word $\alpha$ over the integers let $c_\alpha$ be the curve over the primitive $P_\alpha$ constructed in section 4. I.e.

$$c_\alpha = \begin{cases} d_\alpha & \text{if } g(\alpha) = 1 \\ \text{the curve (4.28) if } g(\alpha) > 1 \end{cases} \qquad (A.2.1)$$

For each finite support function $f: LYN \to \mathbf{N} \cup \{0\}$ let $P^f$ be the ordered product

$$P^f = \prod_{\substack{\alpha \in LYN \\ f(\alpha) > 0}} c_\alpha(f(\alpha)) \qquad (A.2.2)$$

where the wll-ordering is used on $LYN$. It is an easy exercise to show that

$$\mu(P^f) = \sum_{f_1 + f_2 = f} P^{f_1} \otimes P^{f_2} \qquad (A.2.3)$$

A.2.4. Lemma. The $P^f$ form a basis of *NSymm* (as an Abelian group).

Proof (following [12]). First independence. The idea is that by applying the comultiplication a suitable number of times the $c_\alpha(i), i > 1$ finally break up into tensor products of the $c_\alpha(1) = P_\alpha$ (in the 'middle terms') and then independence of the $P_\alpha$ finishes the job.

Define $v_1 = id - e\varepsilon: NSymm \to NSymm$, where $e: \mathbf{Z} \to NSymm$ is the unit morphism and $\varepsilon: NSymm \to \mathbf{Z}$ is the counit. Note that $\nu_1$ kills the constants of *NSymm* and is the identity on the weight $>0$ part. Let

$$\nu_2 = (\nu_1 \otimes \nu_1)\mu: NSymm \to NSymm^{\otimes 2}$$

(comultiplication without the 'trivial parts' $1 \otimes ?$ and $? \otimes 1$) and recursively

$$\nu_n = (\nu_1^{\otimes(n-2)} \otimes \nu_2)\nu_{n-1}: NSymm \to NSymm^{\otimes n}.$$

For instance, if $c(1), c(2), c(3)$ are the first, second, and third term of a $c_\alpha$, then $\nu_2(c(3)) = c(1) \otimes c(2) + c(2) \otimes c(1)$, $\nu_2(c(1)) = 0$, $\nu_3(c(3)) = c(1) \otimes c(1) \otimes c(1)$, $\nu_3(c(2)) = 0$.

Let $S_n$ be the group of permutations on $n$ letters acting on $NSymm^{\otimes n}$ by permuting the factors. For a function of finite support $f: LYN \to \mathbf{N} \cup \{0\}$, let $|f| = \sum_\alpha f(\alpha)$ and let

$$S_f = \prod_{f(\alpha)>0} S_{f(\alpha)} \subset S_{|f|}$$

be the corresponding Young subgroup. Now suppose that the $P^f$ with $|f| < n$, $n \geq 2$, have already been proved to be independent and let $\sum_{|f|\leq n} a_f P^f = 0$. Then

$$0 = \nu_n(\sum_{|f|\leq n} a_f P^f) = \sum_{|f|=n} a_f \sum_{\sigma \in S_n/S_f} \sigma(\bigotimes_\alpha P_\alpha^{\otimes f(\alpha)})$$

so that all the $a_f$ with $|f| = n$ are zero. By the induction hypothesis all the other $a_f$ are then also zero, proving independence of the $P^f$ for all $|f| < n+1$.

Generation is proved by induction on weight. The cases weight zero and weight 1 being obvious. So suppose that $NSymm_k$ is in the linear span of the $P^f$ for all $k < n$ and let $x \in NSymm_n$. Then

$$\mu(x) = 1 \otimes x + x \otimes 1 + \sum_{f,g \neq 0} b_{f,g} P^f \otimes P^g \quad (A.2.5)$$

For certain unique coefficients $b_{f,g}$. (Existence because $\text{wt}(\mu(x)) = \text{wt}(x)$ and the induction hypothesis; uniqueness because of independence). Cocommutativity now ensures

$$b_{f,g} = b_{g,f} \quad (A.2.6)$$

Writing out coassociativity (using the expression (A.2.5) and using (A.2.3)), gives

$$b_{f+g,h} = b_{f,g+h}$$

If follows that there are unique coefficients $c_h$ such that

$$b_{f,g} = c_{f+g}$$

and this in turn gives that $x - \sum_h c_h P^h$ is a primitive (of weight $n$), showing that $x$ is in the linear span of the $P^f$. This proves the lemma.

Now consider the graded cofree cocommutative coalgebra on *LYN* and consider the following morphism of Abelian groups

$$\text{CCoF}(LYN) \longrightarrow NSymm, \quad f \mapsto P^f \tag{A.2.7}$$

By (A.1.5) and (A.2.3) this is a morphism of coalgebras. By Lemma A.2.4 it is surjective. Finally the $\beta_n = \#LYN_n$ satisfy the power series relation (see e.g. [11])

$$(1-t)\prod_{n=1}(1-t^n)^{\beta_n} = 1 - 2t$$

and it follows from the statements above just below (A.1.3) that the rank of $\text{CCoF}(LYN)_n$ is the coefficient of $t^n$ in $(1-t)(1-2t)^{-1}$ which is $2^{n-1} = \text{rank}(NSymm_n)$. Thus (A.2.7) is a homogeneous surjective morphism of graded coalgebras over the integers between graded coalgebras whose free weight $n$ components have the same finite rank. That makes (A.2.7) an isomorphism of coalgebras and hence the graded dual of *NSymm* a free commutative algebra over the integers.